\documentclass[12pt]{article}
\usepackage[utf8]{inputenc}
\usepackage[OT2,T1]{fontenc}
\usepackage{lmodern}
\usepackage[english]{babel}
\usepackage[affil-it]{authblk}
\usepackage{amsmath}
\usepackage{amsfonts,amssymb}
\usepackage{amsthm} 
\usepackage{mathtools}
\usepackage[overload]{empheq}
\usepackage{csquotes}
\usepackage{cite}
\usepackage{wrapfig}
\usepackage[justification=justified, figurename=FIG., format=plain]{caption}
\usepackage{textcomp} 
\DeclareSymbolFont{cyrletters}{OT2}{wncyr}{m}{n}
\DeclareMathSymbol{\Sha}{\mathalpha}{cyrletters}{"58}

\theoremstyle{plain}
\newtheorem{theorem}{Theorem}
\newtheorem{proposition}[theorem]{Proposition}
\newtheorem*{conjecture}{Conjecture}

\theoremstyle{definition}
\newtheorem*{definition}{Definition}

\theoremstyle{remark}
\newtheorem*{remark}{Remark}

\let\OLDthebibliography\thebibliography
\renewcommand\thebibliography[1]{
  \OLDthebibliography{#1}
  \setlength{\parskip}{3.0pt}
  \setlength{\itemsep}{0.1pt plus 0.6ex}
}

\makeatletter
\def\th@plain{
  \thm@notefont{}
  \itshape 
}
\def\th@definition{
  \thm@notefont{}
  \normalfont 
}
\makeatother

\title{\large\bfseries Cyclotomic  
shuffles}
\author[1,2]{Oleg Ogievetsky\thanks{
On leave of absence from P. N. Lebedev Physical Institute, Leninsky Pr. 53,
117924 Moscow, Russia}}
\author[1]{Varvara Petrova}
\affil[\circ]{Aix Marseille Univ, Universit\'e de Toulon, CNRS, CPT, Marseille, France}

\affil[2]{Kazan Federal University, Kremlevskaya 17, Kazan 420008, Russia}
\date{}

\begin{document}
\maketitle
\vspace{-1.2cm}
\begin{abstract}
Analogues of 1-shuffle elements for complex reflection groups of type $G(m,1,n)$ are introduced. A geometric interpretation for $G(m,1,n)$ in terms of rotational permutations of polygonal cards is given. We compute the eigenvalues, and their multiplicities, of the 1-shuffle element in the algebra of the group $G(m,1,n)$. Considering shuffling as a random walk on the group $G(m,1,n)$, we estimate the rate of convergence to randomness of the corresponding Markov chain. We report on the spectrum of the 1-shuffle analogue in the cyclotomic Hecke algebra $H(m,1,n)$ for $m=2$ and small $n$.

\vspace{3mm}
\end{abstract}

\section*{Introduction}
In 1988, N. Wallach considered an element of the group algebra of the symmetric group $S_n$ which is the sum of cycles $(12\dots i)$ where $i$ ranges from $1$ to $n$. He discovered in \cite{NWallach} that the operator of the left multiplication by this element in the group algebra $\mathbb{Z}S_n$ is diagonalizable with eigenvalues
\begin{equation}\label{eq:intro_eigenval}
0,1,2,\dots, n-2,n.
\end{equation}
The sum of cycles $(12\dots i)$ denoted $\Sha_{1,n-1}$ appears in different circumstances and is called 1-shuffle element. In particular, it describes all possible ways of removing the top card from the deck of $n$ cards and inserting it back in the deck at a random position.

Investigating the repeated top-to-random shuffling as a random walk on $S_n$, P. Diaconis et al. \cite{Diaconis_top_to_random} (see also R. Phatarfod \cite{phatarfod}) found that the multiplicity of the eigenvalue $i$ in (\ref{eq:intro_eigenval}) is equal to the number of permutations in $S_n$ with $i$ fixed points, explaining the absence of $n-1$ in (\ref{eq:intro_eigenval}).\\
The $q$-deformation of the result of N. Wallach for the $q$-analogue $^q \Sha_{1,n-1}$ in the Hecke algebra $H_n(q)$ was proposed by G. Lusztig in \cite{Lusztig}. He established that the spectrum of the operator $L_{\Sha_{1,n-1}}$ of the left multiplication by $\Sha_{1,n-1}$ consists of the $q$-numbers
\begin{equation}\label{eq:intro_q_eigenval}
 q^{j-1}[j]_q:=1+q^2+q^4+\dots +q^{2j-2}, \quad j=0,1,\dots, n-2, n.
\end{equation} 

Later, A. Isaev and O. Ogievetsky considered shuffle elements $\Sha_{p,q}$ in the braid group ring $\mathbb{Z}B_n$. With the help of baxterized elements \cite{Ogievetsky_Isaev}, they constructed additive and multiplicative analogues of $\Sha_{p,q}$ in Hecke and Birman-Murakami-Wenzl algebras. The multiplicities of the eigenvalues in (\ref{eq:intro_q_eigenval}) have been established therein by taking the trace of $L_{\Sha_{1,n-1}}:H_n(q) \rightarrow H_n(q)$, using the fact that the $q$-numbers (\ref{eq:intro_q_eigenval}) are linearly independent over $\mathbb{Z}$ as polynomials in $q$. For generic $q$, the multiplicities turn out to be the same as for the symmetric group.

In the present paper we propose polygonal analogues of cards that we call $m$-cards. We introduce elements ${}^{(m)}\Sha_{1,n-1}$, which realise the analogues of top to random shuffling on $m$-cards. The elements ${}^{(m)}\Sha_{1,n-1}$ belong to the group algebra of complex reflection groups of type $G(m,1,n)$.
We adopt the approach of \cite{Diaconis_top_to_random} (for details see \cite{wallah}) and of \cite{Ogievetsky_Isaev} to compute the spectrum and the multiplicities of the eigenvalues of $L_{{}^{(m)}\Sha_{1,n-1}}$. The obtained result for the multiplicities is expressed in terms of the so-called $m$-derangements numbers. Asymptotic convergence to 
randomness in the shuffling the $m$-cards is briefly analysed. We give a preliminary result on the spectrum of $L_{{}^{(m)}\Sha_{1,n-1}}$ in the cyclotomic Hecke algebra $H(m,1,n)$, which is a deformation of the group algebra of the complex reflection group. 

\section{Complex reflection groups $G(m,1,n)$}
A finite complex reflection group is a finite subgroup of $\text{GL}_n(\mathbb{C})$ generated by complex reflections, that is, 
elements $\tau \in \text{GL}_n(\mathbb{C})$ of finite order such that $\text{Ker}(\tau- id)$ is a hyperplane.
The finite complex reflection groups have been classified by Shephard and Todd (1954) into an infinite family of groups $G(m,p,n)$ where $m,p,n$ are positive integers such that $p$ divides $m$, and 34 exceptional groups. 
The group $G(m,1,n)\subset \text{GL}_n(\mathbb{C})$ is formed by all
matrices with exactly one nonzero entry in each row and column; nonzero entries are $m$-th roots of unity; 
it is generated by the elements $t$, $s_1$,$\dots$, $s_{n-1}$ with the defining relations: 
\begin{equation}\label{def2a}
\left\{\begin{array}{ll}
s_{i+1}s_is_{i+1}=s_is_{i+1}s_i & \textrm{for $i=1,\dots,n-2$\ ,}\\[.2em]
s_is_j=s_js_i & \textrm{for $i,j=1,\dots,n-1$ such that $|i-j|>1$}\ ,\\[.2em]
s_i^2=1 & \textrm{for $i=1,\dots,n-1$}
 \end{array}\right.
\end{equation}
and
\begin{equation}\label{def2b}
\left\{\begin{array}{ll}
s_1ts_1t=ts_1ts_1\ ,\\[.2em]
s_it=ts_i & \textrm{for $i>1$\ ,}\\[.2em]
t^m=1\ , & 
\end{array}\right.
\end{equation}
or equivalently by $t_1, \dots, t_n, s_1, \dots, s_{n-1}$ with the defining relations (\ref{def2a}) and
\begin{equation}\label{rel-wr} 
\left\{\begin{array}{ll}
t_i^m=1 & \textrm{for $i=1,\dots,n$\ ,}\\[.2em]
t_i t_j=t_j t_i & \textrm{for $i,j=1,\dots,n$\ ,}\\[.2em]
s_it_j=t_{\pi_i(j)}s_i \ \ & \textrm{for $i=1,\dots,n-1\,,\ j=1,\dots,n\ .$}
\end{array}\right.
\end{equation}
Here $t_1:=t$, $t_{i+1}:=s_it_is_i$ for $i=1,\dots,n-1$ and $\pi_i$ is the transposition $(i,i+1)$. 

In particular, the group $G(1,1,n)$ is isomorphic to the symmetric group $S_n$ and 
the group $G(2,1,n)$ to the hyperoctahedral group $B_n$, the Coxeter group of type $B$.

\vskip .2cm
The group $G(m,1,n)$ also admits the following description, see e.g. 
\cite{Oleg_Induced_repr}.
Let $C_m$ be the cyclic group of order $m$ with a generator $\gamma$. 
Let $E$ be the set $\{ (x,\rho) \; \vert \; x \in \mathbf{n}, \;\rho \in C_m\}$ where $\mathbf{n}=\{1,2,\dots, n\}$.
Define the action of $C_m$ on $E$ by:
\[ \rho \cdot (x,\rho')=(x, \rho \rho'), \quad  \forall \; \rho, \rho' \in C_m \ \text{and}\ \; \forall \; x\in \mathbf{n} .
\]
Let $\text{\sf Perm}(E)$ be the group of all 
permutations of the set $E$. Denote by $\text{\sf Perm}^0(E)$ the subgroup of $\text{\sf Perm}(E)$ consisting of elements $\pi \in \text{\sf Perm}(E)$ such that:
\[ \pi(x,\rho)=\rho \cdot \pi(x,1), \quad  \forall \; \rho \in C_m \ \text{and}\ \; \forall \; x\in \mathbf{n}.
\]
The group $G(m,1,n)$ is isomorphic to $\text{\sf Perm}^0(E)$ via the map $\phi$ such that:
\begin{equation*}
\begin{array}{ll}
\phi(t)(x,1)=\left\{\begin{array}{ll}(x,\gamma) & \textrm{if $x=1$\ ,}\\[.1em]
(x,1) & \textrm{otherwise\, ;}\end{array}\right.  & \\ [1.5em]
\phi(s_i)(x,1)=(\pi_i(x),1)\quad\textrm{for $x\in \mathbf{n}$ and $i= 1,\dots, n-1.$\ }\\
\end{array}
\end{equation*}

We propose a definition of {\it polygonal} cards (or $m$-cards) and interpret $\text{\sf Perm}^0(E)$, and thus $G(m,1,n)$, as the group of \enquote{rotational} permutations of polygonal cards. 

\begin{definition}
An $n$-\emph{deck $\mathfrak{X}$ of $m$-cards} is an $n$-tuple 
$\mathfrak{X}=\left( (x_1, \gamma^{k_1} ),(x_2,\gamma^{k_2}  ), \dots, (x_n, \gamma^{k_n} )\right)$  
where $(x_1,\dots, x_n)$ is a permutation of $(1,\dots, n)$ and $\gamma^{k_j} 
\in C_m$, $ j\in \mathbf{n}$. 
\end{definition}

\noindent
Let $M$ be an $m$-gon with a distinguished vertex. We interpret $(x,1)$ as an $m$-gon $M$, on which the number $x$ is printed, with the distinguished vertex pointing to the north. 
The $m$-card $(x,\gamma^k)$ is the card $(x,1)$ rotated clockwise by an angle $\tfrac{2 \pi}{m}k$ in the plane of $M$.
\setlength{\columnsep}{23pt}
\begin{wrapfigure}[9]{r}{0.57\textwidth}
\vspace{-20pt}
  \begin{center}
   \includegraphics[width=0.57\textwidth]{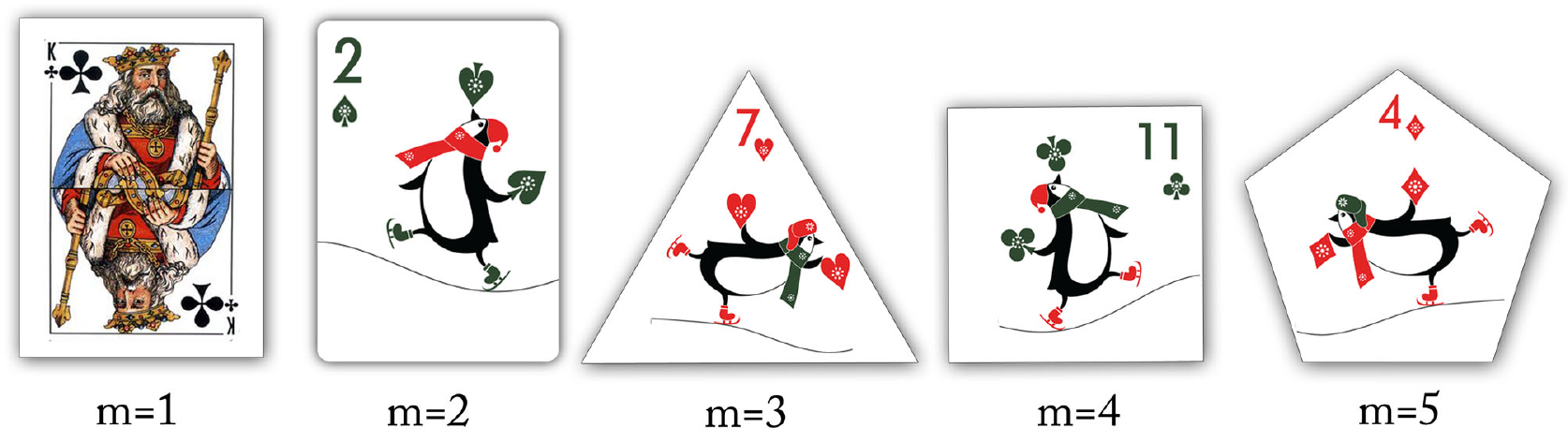}
     \end{center}
  \vspace{-5pt}   
   \caption{m-cards}\label{fig:mcards}
\end{wrapfigure}
\noindent
Abstractly, a 1-card (a usual card) is a point and a 2-card is a segment with two vertices. In the context of card shuffling, we shall slightly transgress this picture and conceive of 1 and 2-cards as the rectangular cards, as illustrated in figure \ref{fig:mcards}.

Let $\left((1,1), \dots , (n,1)\right)$ be an initial deck. We identify a deck $\mathfrak{X}$  
with a \enquote{rotational} permutation of the form: 
$g=\begin{pmatrix} 
(1,1) & \dots & (n,1) \\
(x_1,\gamma^{k_1} ) & \dots & (x_n,\gamma^{k_n} ) 
\end{pmatrix}
$, 
meaning that the $j$-th position in the deck gets occupied by the $x_j$-th card rotated clockwise by an angle $\tfrac{2\pi k_j}{m}$.
The group $\text{\sf Perm}^0(E)$ describes $G(m,1,n)$ as the group of all possible rotational permutations of the $n$-deck of $m$-cards.
In particular, for $i=1,\dots, n-1$, the generator $s_i$ permutes the $i$-th and $(i+1)$-st $m$-cards without any turn.
For $i=1,\dots, n$, the element $t_i$ rotates the $i$-th $m$-card clockwise by an angle $\tfrac{2\pi}{m}$.

\subsection*{Schreier coset graph for the chain of groups $G(m,1,n)$}\label{par:Schreier_coset_graph}
Given a group $G$ with a finite generating set $\mathcal{G}$, and a subgroup $W\subset G$ of finite index,  
the (right) Schreier coset graph of $(G, W;\mathcal{G})$ is an edge-labeled graph whose vertices are the right $W$-cosets  
and edges are of the form $(Wg, Wgx)$ where 
$x\in \mathcal{G}$ and $g \in G$; the edge $(Wg, Wgx)$ has label $x$. 
The graph depicts the action of $G$ on cosets of $W$ by right multiplication and is equivalent to the coset table obtained by the Todd-Coxeter algorithm \cite{coset_enum}. 

\begin{figure}[htp]
\centering
\includegraphics[scale=0.6]{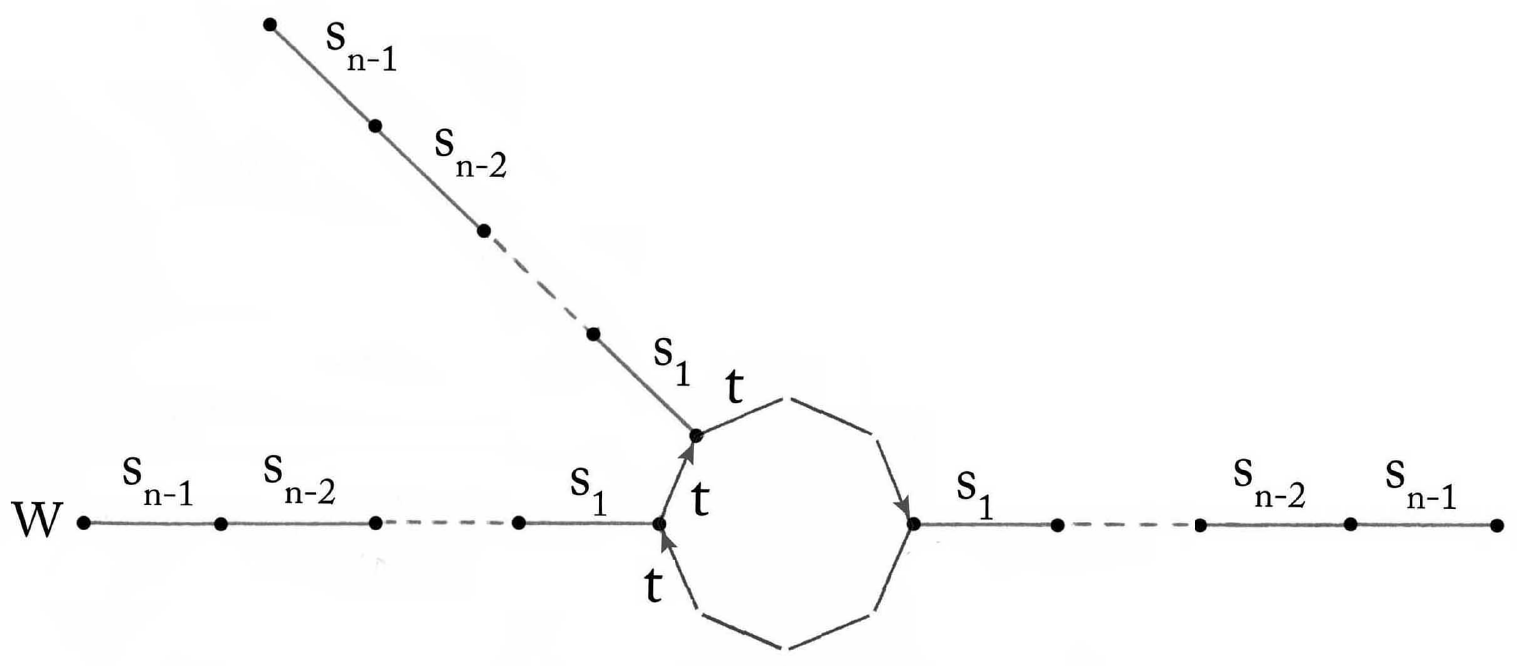}
\caption{Schreier coset graph for $G(m,1,n)$ with respect to $W$ \cite{Oleg_Induced_repr}.}\label{fig:shr_gr_G}
\end{figure}
Figure \ref{fig:shr_gr_G} is the Schreier coset graph 
of $\left(G(m,1,n),W;\mathcal{G}\right)$ where $\mathcal{G}=\{ t, s_1, \dots, s_{n-1}\}$ and $W$ is generated by the elements $t, s_1, \dots, s_{n-2}$. An unoriented edge represents a pair of edges with opposite directions.  For simplicity, only non-trivial actions are drawn (if $Wgx=Wg$ then $(Wg, Wgx)$ is a loop and we do not draw it). At each vertex of the oriented $m$-gon in the middle of the figure starts a tail with $n-1$ edges. All tails are identically (edges-)labelled.\\
The graph provides a \enquote{normal form with respect~to~$W$} \cite{Oleg_Induced_repr}. 
\begin{proposition}\label{prop:normalform}
Any element $x \in G(m,1,n)$ can be uniquely written in the form:
$ x=w s_{n-1} \cdots s_1 t^{\alpha}s_1\cdots s_j
$ with $j \in \{0,\dots, n-1\}$, $\alpha\in\{0, \dots, m-1 \}$ and ${w\in W\simeq G(m,1,n-1)}$ (by convention, the empty product is equal to 1).
\end{proposition}
Starting with a normal form (powers of $t$) for $G(m,1,1)\simeq C_m$, we build recursively \cite{Oleg_Induced_repr} the global normal form  for elements of $G(m,1,n)$ for any $n$. We have an ascending chain:
\[G(m,1,0) \subset G(m,1,1) \subset G(m,1,2) \subset \dots
\]

In the computational context, in order to build an algorithm that takes as input a word in generators of ${G(m,1,n+1)}$ and their inverses, and returns its normal form one only needs to know the normal form of the words $\psi_{k j}^{(n)}x=s_{n} \cdots s_1 t^{k}s_1\cdots s_j x$ for all generators $x$ of $G(m,1,n+1)$, $k=0,\dots, m-1$ and $j=0, \dots, n.$ 
We established the following recursive list of rewriting rules (the equalities are to be understood as \enquote{replace the left hand side by the right hand side}). 
\begin{proposition}
The rewriting rules for $G(m,1,n+1)$ with $n>1, m \geqslant 1$ include (additionally to the defining relations (\ref{def2a})-(\ref{def2b}))
the rewriting rules for $G(m,1,n)$ and the following:
\begin{flalign*}[left=\empheqlbrace \;]
&s_n s_{n-1} \cdots s_j s_n= s_{n-1}s_n s_{n-1} \cdots s_j, && j=1,...,n-2, \\
&\psi_{k j}^{(n)} s_n = s_{n-1}\psi_{k j}^{(n)}, && k=0,...,m-1, \; j=0,1,..., n-2.
\end{flalign*}
For $G(m,1,2)$ the rewriting rules (in addition to the defining relations) are: 
\[
s_1 t^k s_1 t=t s_1 t^k s_1, \quad \forall \; k=2,...,m-1.
\]
\end{proposition}

\section{Shuffle elements in $G(m,1,n)$ group algebra}
We start by recalling that a $(p,q)$-\emph{shuffle}, where $p$ and $q$ are non-negative integers, is a permutation ${\sigma\in S_{p+q}}$ of the set $\{1,\dots, p+q\}$ such that $\sigma(1)<\sigma(2) < \dots < \sigma(p)$ and $\sigma(p+1)<\dots < \sigma(p+q)$. 
When $p=1$ or $q=1$, the $(p,q)$-shuffle is called 1-shuffle.
The sum over all $(p,q)$-shuffles with $p,q$ fixed is called shuffle element and is denoted by~$\Sha_{p,q}$. 

Let $w^{\uparrow \ell}$ denote the image of $w \in S_n$ under the homomorphism $S_n \to S_{n+\ell}$ sending $\pi_i=(i, i+1)$ to $\pi_{i+\ell}$, $i=1,\dots, n-1$.  
Shuffles $\Sha_{p,q}$ can be defined inductively by any of the recurrence relations (analogues of the Pascal rule): 
\begin{eqnarray*}
\Sha_{p,q}&=&\Sha_{p-1,q}+ \Sha_{p, q-1} \pi_{p+q-1} \cdots \pi_q\\
&=&\Sha_{p,q-1}^{\uparrow 1}+ \Sha_{p-1,q}^{\uparrow 1} \pi_1 \cdots \pi_q\ ,
\end{eqnarray*}
with $\Sha_{0,q}=\Sha_{q,0}=1$ for any $q \geqslant 0$. 
Let $\tilde{\Sigma}_n:=\sum_{\sigma \in S_n}\sigma$ be the symmetrizer, $\sigma \tilde{\Sigma}_n = \tilde{\Sigma}_n,  \; \forall \sigma \in S_n$.  
Then 
$\tilde{\Sigma}_{q+p}=\Sha_{p,q} \tilde{\Sigma}_q \tilde{\Sigma}_p^{\uparrow q}$, in particular,
\begin{equation}\label{eq:rec_sym_1sh}
\tilde{\Sigma}_n=\Sha_{1,n-1} \tilde{\Sigma}_{n-1}.
\end{equation} 

By their combinatorial structure and relationship with the structure of products, shuffle elements appear in many constructions in homotopy theory and higher category theory.
In particular, they are involved in the product of simplices \cite{Mac_Lane_def_shuffle}, in the study of multiple polylogarithms 
\cite{Bowman, Borwein}, 
in the bialgebras of type 1 \cite{N}, 
in the construction of the Hopf algebra structure on the tensor algebra $T(V)$ of a vector space $V$ \cite{GO}, 
in the construction of the standard complex for quantum Lie algebras \cite{Ogiev_Isaev2004brst}, etc.

\vspace{1mm}
In the present paper we consider the analogues of the 1-shuffle elements $\Sha_{1,n-1}$ in $\mathbb{Z}G(m,1,n)$ of the form
\[{}^{(m)}\Sha_{1,n-1}:=\Sha_{1,n-1} \cdot {}^{(m)}\mathcal{T}_n,
\]
where $\Sha_{1,n-1}= 1+s_{n-1}+s_{n-2}s_{n-1}+...+s_1 \cdots s_{n-1}$ and
${}^{(m)}\mathcal{T}_n:=1+t_n+t_n^2+...+t_n^{m-1}.$ 
The element ${}^{(m)}\Sha_{1,n-1}$ describes all possible ways of rotating the $n$-th $m$-card and inserting it at random position in the deck; we found that ${}^{(m)}\Sha_{1,n-1}$ satisfies:
\begin{equation}\label{eq:car_id_G}
{}^{(m)}\Sha_{1,n-1}^2 = {}^{(m)}\Sha_{1,n-1}\left(m + {}^{(m)}\Sha_{1,n-2}\right).
\end{equation}
Let ${}^{(m)}\tilde{\Sigma}_n$ denote the sum of all elements of the group $G(m,1,n)$. The recurrent relation similar to (\ref{eq:rec_sym_1sh}) holds:
\begin{equation*}
{}^{(m)}\tilde{\Sigma}_n = {}^{(m)}\Sha_{1,n-1} \cdot {}^{(m)}\tilde{\Sigma}_{n-1}.
\end{equation*}
Making use of the relation (\ref{eq:car_id_G}) we prove by induction on $n$ that
\begin{equation}\label{eq:prod_sh_eq_symm}
\prod_{i=0}^{n -1} \left({}^{(m)}\Sha_{1,n-1} - i m \right)  = {}^{(m)}\tilde{\Sigma}_n\ . 
\end{equation}
Since ${}^{(m)}\Sha_{1,n-1} \cdot {}^{(m)}\tilde{\Sigma}_n  = mn\cdot {}^{(m)}\tilde{\Sigma}_n$, we obtain, multiplying (\ref{eq:prod_sh_eq_symm})  
by  $\left({}^{(m)}\Sha_{1,n-1} - nm  \right)$, 
\begin{equation}\label{eq:min_pol}
\prod_{i=0}^{n}\left({}^{(m)}\Sha_{1,n-1}-im\right)=0.
\end{equation}
For a given $x \in \mathbb{Z}G(m,1,n)$, let $L_x$ be the map $\mathbb{Z}G(m,1,n) \rightarrow \mathbb{Z}G(m,1,n)$, $y \mapsto xy$. 
We calculate the multiplicities of the eigenvalues of $L_{{}^{(m)}\Sha_{1,n-1}}$ following the approach Garsia and Wallach developed for the symmetric group in \cite{wallah}. 
Let $\alpha$ and $\beta$ be two words in some alphabet. We denote by $\alpha\doublecup\beta$ the formal sum of all words that can be obtained by shuffling 
$\alpha$ and $\beta$ as it is done with two decks of cards. 
We define for $0\leqslant a\leqslant n$: 
$$B_a = \sum_{\alpha\in G(m,1,a)}\alpha\doublecup\beta_{a,n}\ ,$$
where $\alpha$ is viewed as an $a$-deck of $m$-cards and $\beta_{a,n}$ is the deck $((a + 1, 1), (a + 2, 1), \dots , (n, 1))$.
In particular, $B_0 = 1$. Since $B_1={}^{(m)}\Sha_{n-1,1}=w\cdot {{}^{(m)}\Sha_{1,n-1}}\cdot w^{-1}$ where $w$ is the longest element of 
the symmetric group $S_n=G(1,1,n)$, it has the same spectrum as ${}^{(m)}\Sha_{n-1,1}$. The element $B_a$ is composed of $\frac{n!}{(n-a)!}m^a$ decks that may be grouped according to the number and the angle of the top card. We get $am+1$ groups labelled by $(1,\rho),(2,\rho),\dots,(a,\rho)$ and $(a+1,1)$, with $\rho\in C_m$.
Upon multiplication by $B_1$, each group $(k,\rho)$, $1\leqslant k\leqslant a$, $\rho\in C_m$, yields a term $B_a$, and the group $(a+1,1)$ gives the term $B_{a+1}$. 
Thus, $B_1B_a=amB_a+B_{a+1}$ or $B_{a+1}=(B_1-am)B_a$, $0\leqslant a\leqslant n-1$. Iterating, we obtain 
\begin{equation}\label{itefoba}
B_a=B_1(B_1-m)(B_1-2m)\cdots (B_1-(a-1)m)\ ,\ 0\leqslant a\leqslant n\ .\end{equation}
Applying the Stirling inversion to (\ref{itefoba}), we find, for $0\leqslant k \leqslant n$: 
\begin{equation}\label{eq:B1k_sum_proj}
B_1^k   =m^k \sum_{i=0}^n i^k \mathbf{E}_i\ ,
\ \ \text{where}\ \ 
\mathbf{E}_i=\frac{1}{i!} \sum_{a=i}^n \frac{(-1)^{a-i}}{(a-i)!m^a} \, B_a\ .  
\end{equation}
Substituting into the identity $B_1B_1^k=B_1^{k+1}$ the expressions (\ref{eq:B1k_sum_proj}) for $B_1^k$ and $B_1^{k+1}$ we obtain:
\begin{equation*}
\sum_{i=0}^n i^k\left( B_1 \mathbf{E}_i - im \mathbf{E}_i \right) =0\ , \quad k=0,\dots, n\ .
\end{equation*}
Since the matrix $M_{ik}:=\{i^k\}_{i,k=0}^n$ is invertible, we get
$B_1 \mathbf{E}_i = im \mathbf{E}_i$, $i=0 ,\dots, n$. 

Let $\mathcal{B}$ be the $\mathbb{Q}$-subalgebra of $\mathbb{Q}G(m,1,n)$, generated by $B_1$. 
By (\ref{itefoba}), $\mathcal{B}$ is a linear span of $B_a$, 
$a=0, \dots, n$. By (\ref{eq:B1k_sum_proj}),  
$\mathcal{B}$ is a linear span of $\mathbf{E}_a, a=0, \dots, n$.
Therefore $\mathbf{E}_i$ and $\mathbf{E}_j$ commute. We have
\[  B_1\mathbf{E}_i \mathbf{E}_j = mi \mathbf{E}_i \mathbf{E}_j = mj \mathbf{E}_i \mathbf{E}_j\ , \quad 0 \leqslant i,j \leqslant n\ ,
\]
whence $\mathbf{E}_i \mathbf{E}_j = 0 , \; 0 \leqslant i\neq j \leqslant n$.
On the other hand, the equality $\sum_{i=0}^n \mathbf{E}_i=1$ obtained from (\ref{eq:B1k_sum_proj}) for $k=0$ implies that $\mathbf{E}_i^2 = \mathbf{E}_i, \; i=0, \dots, n$.
Operators $\mathbf{E}_i$, $i=0,\dots, n,$ are therefore the orthogonal projections onto the eigenspaces of $B_1$. 
The multiplicity of the eigenvalue $im$ in (\ref{eq:min_pol}) is thus given by the trace of the matrix $L_{\mathbf{E}_i}$. Since the trace of any element except the identity vanishes in the regular representation, we get:
\begin{equation*}
tr(L_{\mathbf{E}_i})=\frac{1}{i!} \sum_{a=i}^n \frac{(-1)^{a-i}}{(a-i)!m^a} tr(L_{B_a})
=\binom{n}{i}m^{n-i}(n-i)!\sum_{a=0}^{n-i}\frac{(-1)^a}{m^a a!}\ .
\end{equation*}
The obtained multiplicities have the combinatorial meaning. 
We shall say that an element $g \in G(m,1,n)$ has $i$ fixed points if it fixes $i$ $m$-cards.
For a given $n$ the $m$-\emph{derangement number} is the number $d_{m,n}$ of elements in $G(m,1,n)$ with no fixed points. We have
\begin{equation*}
d_{m,n}=m^n n! \sum_{k=0}^{n}\frac{(-1)^k}{k!m^k}\ .
\end{equation*}
We summarize the obtained results.
\begin{theorem} The operator $L_{{}^{(m)}\Sha_{1,n-1}}: \mathbb{Z}G(m,1,n) \rightarrow \mathbb{Z}G(m,1,n)$ is diagonalizable with eigenvalues $im$, 
$i=0,1,\dots,n-1,n$. The multiplicity of the eigenvalue $im$ equals $\binom{n}{i}d_{m,n-i}$, the number of elements of $G(m,1,n)$ with $i$ fixed points. 
\end{theorem}
\begin{remark}
The spectrum of $L_{{}^{(1)}\Sha_{1,n-1}}$ is $\{0,1,\dots, n-2, n\}$, because $d_{1,1}=0$; the spectrum has a gap at $n-1$. For $m>1$, the gap disappears. This result proves important in estimating how close $k$ repeated shuffles, randomly picked among the set of elements of $G(m,1,n)$ whose sum is ${}^{(m)}\Sha_{1,n-1}$, get the $n$-deck of $m$-cards to being randomized. Indeed, under a suitable normalization, the matrix of $L_{{}^{(m)}\Sha_{1,n-1}}$ is a matrix $P$ of transition probabilities. As it belongs to the convex envelope of the set of permutation matrices, the matrix $P$ is doubly stochastic. The underlying Markov chain is ergodic and its unique stationary probability distribution is uniform.
Using the spectral representation formula for the $k$-step transition matrix $P^k$ \cite{Seneta, spec_dec}, one can show that the convergence to the stationarity occurs exponentially fast with the rate of convergence governed by the second largest eigenvalue: there exists a constant $\alpha>0$ such that:
\[ \forall\; k \geqslant 1, \quad \left| P^k_{ij}-\frac{1}{m^n n!} \right| \leqslant 
\left\{
\begin{array}{ccl}
\alpha \left( \frac{n-2}{n} \right)^k & \mbox{if} & m=1\ ,\\
\alpha \left( \frac{n-1}{n} \right)^k & \mbox{if} & m>1\ .\\
\end{array}
\right.
\]
\end{remark}

\paragraph{Cyclotomic algebras.}
The cyclotomic Hecke algebra $H(m,1,n)$ is the associative algebra over $\mathbb{C}\left[q^{\pm}, v_1^{\pm}, \dots, v_m^{\pm} \right]$ generated by $\tau,\sigma_1,\dots,\sigma_{n-1}$ with the defining relations:
\begin{flalign*}[left=\empheqlbrace  \;]
&\sigma_i\sigma_{i+1}\sigma_i=\sigma_{i+1}\sigma_i\sigma_{i+1} &\hspace{1cm}& \textrm{for all $i=1,\dots,n-2 \ ,$}\\
&\sigma_i\sigma_j=\sigma_j\sigma_i && \textrm{for all $i,j=1,\dots,n-1$ such that $|i-j|>1 \ ,$}\\
&\tau\sigma_1\tau\sigma_1=\sigma_1\tau\sigma_1\tau\, &&\\
&\tau\sigma_i=\sigma_i\tau && \textrm{for $i>1 \ ,$}\\
 & \sigma_i^2=(q-q^{-1})\sigma_i+1 & &  \textrm{for all $i=1,\dots,n-1 \ ,$}\\
 & (\tau-v_1)\cdots(\tau-v_m)=0. & & 
\end{flalign*}
The algebra $H(m,1,n)$ is a deformation of the complex reflection group algebra $\mathbb{C} G(m,1,n)$. 
The representation theory of $H(m,1,n)$ has been constructed in \cite{ref1AK}; the representation theory of $G(m,1,n)$ and $H(m,1,n)$ was reobtained in \cite{Ogievetsky_Poulain_JM_elts,Ogievetsky_Poulain_JM_Gm1n} from the analysis of the spectrum of the so called Jucys-Murphy elements, following the approach of \cite{OV}. 

Consider the algebra $H(2,1,n)$ with parameters $v_1=p$, $v_2=-p^{-1}$. Let elements $\tau_i$ be recursively defined by, cf (\ref{rel-wr}): 
\begin{equation*}
\tau_i=\sigma_{i-1}^{-1} \tau_{i-1} \sigma_{i-1}, \quad i=2,\dots, n \ ,
\end{equation*}
with $\tau_1=\tau$.
We propose the following analogue of ${}^{(2)}\Sha_{1,n-1}$ for the algebra $H(2,1,n)$:
\begin{equation}\label{eq:shHecke2}
{}^{(2)}\Sha_{1,n-1;q,p}:=\left(1+q \sigma_{n-1}+q^2\sigma_{n-2}\sigma_{n-1}+ \dots + q^{n-1} \sigma_{1}\cdots \sigma_{n-1}\right) \left( 1+p \tau_n \right).
\end{equation}
Our calculations of the spectrum of (\ref{eq:shHecke2}) in various representations of $H(2,1,n)$ for small $n$ lead to the following conjecture. 

\begin{conjecture}
The spectrum of the operator of the left multiplication by ${}^{(2)}\Sha_{1,n-1;q,p}$ in the algebra $H(2,1,n)$ consists of the values:
\begin{equation*}
\left(p^2+ q^{2(j-\ell)} \right) q^{\ell-1}[\ell]_q \ ,
\end{equation*}
where $q^{\ell-1}[\ell]_q$ are the $q$-numbers defined in (\ref{eq:intro_q_eigenval}), $\ell=0,1,\dots, n$ and  $j=\ell, {\ell+1}, \dots, n-2, n$.
\end{conjecture}
\begin{remark}
We numerically observed (for $n\leqslant 5$) that the eigenvalues 
$0$, $(p^2+q^2)q^{n-2}[n-1]_q$ and $(p^2+1)q^{n-1}[n]_q$ have the same multiplicities as their respective counterparts $0, 2(n-1), 2n$ in the classical situation, 
for the group $G(m,1,n)$.
This is not the case for the other eigenvalues.
\end{remark} 

\vskip .2cm\noindent {\bf Acknowledgements.}
The work of O. O. was supported by the Program of Competitive Growth of Kazan Federal University and by the grant RFBR 17-01-00585.
The work of V.~P. has been carried out thanks to the support of the A*MIDEX grant (\textnumero~ANR-11-IDEX-0001-02) funded by the French Government \enquote{Investissements d’Avenir} program.

\end{document}